\documentclass{amsart}

\usepackage{amssymb}

\usepackage[utf8]{inputenc}
\usepackage[T1]{fontenc}

\usepackage[english,francais]{babel}

\newtheorem{theorem}{Theorem}[section]

\newtheorem{e-proposition}[theorem]{Proposition}

\newtheorem{e-definition}[theorem]{Definition\rm}

\newtheorem{theoreme}{Th\'eor\`eme}[section]

\newtheorem{proposition}[theoreme]{Proposition}

\newtheorem{definition}[theoreme]{D\'efinition\rm}

\begin{document}
\title[Approximations and Bernoulli percolation]{Approximations of standard equivalence relations and Bernoulli percolation at $p_u$}
\author[D. Gaboriau]{Damien Gaboriau}
\address{CNRS, Unit\'e de math\'ematiques Pures et Appliqu\'ees, ENS-Lyon, Universit\'e de Lyon, France}
\email{damien.gaboriau@ens-lyon.fr}
\author[R. Tucker-Drob]{Robin Tucker-Drob}
\address{Department of Mathematics, Rutgers University, Piscataway NJ}
\email{rtuckerd@math.rutgers.edu}
\keywords{ergodic theory, Bernoulli percolation, orbit-equivalence, strong ergodicity}

\maketitle

\begin{abstract}
\selectlanguage{english}
The goal of this note is to announce certain results (to appear in \cite{Gab-Tuck=Approx-in-prepa}) in orbit equivalence theory, especially concerning the approximation of p.m.p. standard equivalence relations by increasing sequence of sub-relations, with applications to the behavior of the Bernoulli percolation on graphs at the threshold $p_u$.

\vskip 0.5\baselineskip

{\selectlanguage{francais} 

\noindent{R\'esum\'e en Fran\c cais.}
 Approximations de relations d'équivalence standards et percolation de Bernoulli à $p_u$. }
Le but de cette note est d'annoncer certains résultats (à paraître dans \cite{Gab-Tuck=Approx-in-prepa}) d'équivalence orbitale, concernant notamment la notion d'approximation par suite croissante de sous-relations, avec applications au comportement en $p_u$ de la percolation de Bernoulli sur les graphes de Cayley.

\end{abstract}

\section{Version fran\c{c}aise abr\'eg\'ee}

La notion de relation d'équivalence standard \emph{hyperfinie} (i.e. réunion croissante de sous-relations standards finies) joue un rôle fondamental en théorie de l'équivalence orbitale. 
Plus généralement, on peut considérer la notion d'\emph{approximation} d'une relation d'équivalence mesurée standard ${\mathcal{R}}$, i.e. la possibilité d'écrire ${\mathcal{R}}$ comme réunion croissante d'une suite de sous-relations d'équivalence standards ${\mathcal{R}}=\bigcup_{n\in \mathbb{N}}\nearrow {\mathcal{R}}_n$. Une telle approximation est \emph{triviale} s'il existe une partie borélienne $A$ non triviale sur laquelle les restrictions coïncident à partir d'un certain rang ${\mathcal{R}}_n{\upharpoonright} A={\mathcal{R}}{\upharpoonright} A$.
Nous établissons des conditions sous lesquelles les approximations de certaines relations d'équivalence sont nécessairement triviales.

\begin{theoreme}\label{th:red herring fr}
Soit $G$ un groupe engendré par deux sous groupes (infinis) de type fini commutant $H$ et $K$.
Considérons une action libre \emph{préservant la mesure de probabilité} (\emph{p.m.p.}) sur l'espace borélien standard 
$G\overset{\alpha}{\curvearrowright}(X,\mu)$ telle que $H$ agit de manière fortement ergodique et $K$ de manière ergodique.
Alors la relation d'équivalence engendrée ${\mathcal{R}}_{\alpha}$ n'admet aucune approximation non triviale.
\end{theoreme}

Puisque les actions Bernoulli des groupes non moyennables sont automatiquement fortement ergodiques, ce résultat a des conséquences en théorie de la percolation de Bernoulli sur les graphes de Cayley. Pour des informations concernant les liens entre équivalence orbitale et percolation, on peut consulter \cite{Gab05}.
En fait, le couplage standard permet de traduire l'étude relative aux variations du paramètre de rétention $p\in [0,1]$ de la percolation en l'étude d'une famille croissante de relations d'équivalence standards p.m.p. $({\mathcal{R}}_p)_{p\in [0,1]}$ telle que pour tout $q\in [0,1]$, on a  ${\mathcal{R}}_q=\bigcup_{p<q} \nearrow {\mathcal{R}}_p$.
Le paramètre critique $p_u$ ({\it cf.} \cite{Haggstrom-Peres=1999}) est l'infimum des $p$ pour lesquels on peut trouver une partie borélienne non négligeable $A$ sur laquelle les restrictions ${\mathcal{R}}_1{\upharpoonright} A$ et ${\mathcal{R}}_p{\upharpoonright} A$ coïncident (de tels $p$ sont dits appartenir à la \emph{phase d'unicité}).
Pour les groupes dont les actions Bernoulli n'admettent pas approximation non triviale, alors $p_u$ lui-même n'appartient pas à la phase d'unicité. \emph{C'est le cas des groupes qui apparaissent dans le théorème~\ref{th:red herring fr}}. Des conditions d'exhaustion par des sous-groupes distingués en un sens faible nous permettent d'élargir encore la famille de nouveaux exemples.

Les notions de \emph{dimension géométrique} et de \emph{dimension approximative} d'une relation d'équivalence mesurée ont été introduites dans \cite[section 5]{Gab02}, où il est démontré qu'une non annulation du $d$-ième nombre de Betti $\ell^2$ fournit une minoration par $d$ de ces deux notions de  dimension. La première est analogue à la notion de dimension géométrique pour un groupe et la deuxième est le minimum des $\liminf$ des dimensions géométriques le long des suites approximantes.
Pour les relations non approximables, les deux notions de dimension coïncident manifestement.
On peut alors exhiber des familles de groupes possédant des actions de dimensions approximatives variables.
\selectlanguage{english}
\section*{English version}
\section{Bernoulli bond percolation}
Let ${\mathcal{G}}=(G,\mathtt{E})$ be a Cayley graph for a finitely generated group $G$.
The \emph{Bernoulli bond percolation} on ${\mathcal{G}}$, with retention parameter $p\in [0,1]$, considers the i.i.d. assignment to each edge in $\mathtt{E}$ of the value $1$ (open) with probability $p$ and of the value $0$ (closed) with probability $1-p$. 
The number of infinite \emph{clusters} (connected components of open edges), for the resulting probability measure $\mathbf{P}_{p}$ on $\{0,1\}^{\mathtt{E}}$,  is $\mathbf{P}_p$-a.s. either 0, 1 or $\infty$. 
There two critical values, $0<p_{c}({\mathcal{G}})\leq p_{u}({\mathcal{G}})\leq 1$, depending on the graph which govern the three regimes, as summarized in the following picture (see \cite{Haggstrom-Peres=1999}):
\vspace{-5pt} 
\begin{center}                                                               
\begin{tabular}{lcccccr}
    & \raisebox{-5pt}[0cm][0cm]
    {all finite} & & \hskip-5pt \raisebox{-5pt}[0cm][0cm]{$\infty^{\textrm{ly}}$ many $\infty$ 
    clusters} \hskip-5pt  & & \hskip-5pt  \raisebox{-5pt}[0cm][0cm]{a unique $\infty$ cluster} \hskip-5pt &
    \\
    \multicolumn{7}{l}{\hrulefill}\hskip 2.5pt
    \\
    \hskip -1.5pt \raisebox{3.6ex}[0cm][0cm]{$\vert$} & & 
    \raisebox{3.6ex}[0cm][0cm]{$\vert$} &  &  \raisebox{3.6ex}[0cm][0cm]{$\vert$}
    & &  \raisebox{3.6ex}[0cm][0cm]{$\vert$}\phantom{\hskip 1pt}
    \\
    \hskip -2pt \raisebox{10pt}[0cm][0cm]{$0$} & & \hskip-5pt  \raisebox{10pt}[0cm][0cm]{$p_{c}({\mathcal{G}})$} &  &  \hskip-5pt  
    \raisebox{10pt}[0cm][0cm]{$p_{u}({\mathcal{G}})$} & \raisebox{10pt}[0cm][0cm]{uniqueness phase} & \raisebox{10pt}[0cm][0cm]{$1$}
\end{tabular}
\end{center}
\vskip-12pt
While it is far from being entirely understood, there are some partial results concerning the situation at the threshold $p=p_{u}$ and our Theorem~\ref{th:nonuniqeness at pu} contributes to this study. 

For groups with infinitely many ends, $p_{u}=1$  \cite{LP05}; thus the percolation 
    at $p=p_{u}$ belongs to the uniqueness phase.
At the opposite, the percolation at the threshold $p=p_{u}$ belongs to the nonuniqueness phase (and thus $p_{u}<1$) for all Cayley graphs of infinite groups with Kazhdan's property (T) \cite{Lyons-Schramm=1999}. 
Y.~Peres \cite{Peres-2000=pu-products} proved that for a non-amenable direct product of infinite groups $G=H\times K$, and for a Cayley graph associated with a generating system $S=S_H\cup S_K$ with $S_H\subset H$ and $S_K\subset K$, then the percolation at $p_u({\mathcal{G}})$ doesn't belong to the uniqueness phase. 
We extend this result to a larger family of groups than direct products, and to any of their Cayley graphs.

\begin{theorem}[Non uniqueness at $p_u$]\label{th:nonuniqeness at pu}
Let  $G$ be a non-amenable group generated by two commuting infinite and finitely generated subgroups $H$ and $K$. Then for every Cayley graph ${\mathcal{G}}$ of $G$, the percolation at $p_u({\mathcal{G}})$ doesn't belong to the uniqueness phase. 
\end{theorem}
The same result holds when $G$ admits an infinite normal subgroup $H$ such that the pair $(G,H)$ has relative property (T). This has also been observed by C. Houdayer (personal communication).
Using some weak forms of normality we can extend the scope of our theorem, for instance when $G$ is a nonamenable (generalized) Baumslag-Solitar group (see Theorem~\ref{th: BS-groups}), or a nonamenable HNN-extension of $\mathbb{Z} ^n$ relative to an isomorphism between two finite index subgroups.

Theorem~\ref{th:nonuniqeness at pu} follows from a general result on approximations of standard probability measure preserving equivalence relations (Th.~\ref{th:red herring fr}).
We refer to \cite{Gab05} for general informations concerning connections between equivalence relations  and percolation on graphs, and references therein.

\section{Approximations of standard equivalence relations}

Let ${\mathcal{R}}$ be a standard probability measure preserving (\emph{p.m.p.}) equivalence relation on the atomless probability standard Borel space $(X,\mu)$. See \cite{FM77a} for a general axiomatization of this notion.

\begin{definition}[Approximations]
An \emph{approximation} $({\mathcal{R}}_n)$   to ${\mathcal{R}}$ is an exhausting increasing sequence of sub-equivalence relations: $\bigcup_{n\in \mathbb{N}}\nearrow {\mathcal{R}}_n={\mathcal{R}}$.
An approximation is \emph{trivial} if there is some $n$ and a non-negligeable Borel subset $A\subset X$ on which the restrictions coincide: ${\mathcal{R}}_n {\upharpoonright} A={\mathcal{R}} {\upharpoonright} A$. 
We say that ${\mathcal{R}}$ is \emph{non-approximable} if every approximation is trivial.
An action $G\overset{\alpha}{\curvearrowright}(X,\mu)$ is \emph{approximable} if its orbit equivalence relation ${\mathcal{R}}_{G}:=\{(x,\alpha(g)(x)): x\in X, g\in G\}$ is approximable.
\end{definition}

For instance, all free p.m.p. actions of a non-finitely generated group are approximable.
Finite standard equivalence relations are non-approximable.

\begin{proposition}[Approximable equivalence relations]
The following are examples of approximable equivalence relations.
\begin{enumerate}
\item \label{it: OW}
Every aperiodic p.m.p. action of an (infinite) amenable group is approximable by a sequence of sub-equivalence relations with finite classes. 

\item\label{it:non-strong-erg} Every ergodic non-strongly ergodic p.m.p. equivalence relation admits an approximation by 
${\mathcal{R}}_n$ with diffuse ergodic decompositions. 

\item\label{it:free prod and approx} Any free product ${\mathcal{R}}=\mathcal{A}*\mathcal{B}$ of aperiodic p.m.p. equivalence relations is approximable.
\end{enumerate}
\end{proposition}
Item~(\ref{it: OW}) follows from Ornstein-Weiss theorem \cite{OW80}.
Item~(\ref{it:non-strong-erg}) relies heavily on results of Jones-Schmidt \cite{Jones-Schmidt-1987}.
Recall that \emph{strong ergodicity}, a reinforcement of ergodicity introduced by K.~Schmidt, 
requires that:
for every sequence $(A_n)$ of Borel subsets of $X$ such that the symmetric differences satisfy $\lim_{n\to \infty}\mu(A_n\Delta g.A_n)=0$ for each $g\in G$, we must have $\mu(A_n)(1-\mu(A_n))\to 0$.
Item~(\ref{it:free prod and approx}) will be developed in \cite{Gab-Tuck=Approx-in-prepa}.

\begin{proposition}[Non-approximable equivalence relations]
The following are examples of non-approximable equivalence relations.
\begin{enumerate}
\item Every p.m.p. action of a Kazdhan property (T) group is non-approximable. 
\item Every free p.m.p. actions of $\mathrm{SL}(2,\mathbb{Z})\ltimes \mathbb{Z} ^2$, where $\mathbb{Z}^2$ acts ergodically, is non-approximable. More generally free actions of relative property (T) pairs $(G,H)$ where $H$ is normal, infinite and acts ergodically.

\end{enumerate}
\end{proposition}

We prove the following effective version of Th.~\ref{th:red herring fr}.
\begin{theorem}[Effective non-approximability]\label{thm:first}

Let $G$ be a countable group generated by two commuting subgroups $H$ and $K$. Consider a p.m.p. action $G\curvearrowright (X,\mu )$ of $G$ in which $H$ acts strongly ergodically and $K$ acts ergodically.
Let $\mathcal{E}$ is any Borel subequivalence relation of $\mathcal{R}_G$. For each $g\in G$, set $A_g:= \{ x\in X : g x \, \mathcal{E} \, x \}$.
Let $S$ and $T$ be generating sets for $H$ and $K$ respectively. For every $\epsilon >0$, there exists $\delta >0$ such that if $\mathcal{E}$ satisfies:
\begin{enumerate}
\item[(i)] $\mu (A_s)
> 1-\delta$ for all $s\in S$, and
\item[(ii)] $\mu (A_t)  
> \epsilon$ for all $t\in T$,
\end{enumerate}
then there exists a
Borel set $B\subseteq X$, with $\mu (B)>1-\epsilon$, where the restrictions coincide: $\mathcal{E}{\upharpoonright} B= \mathcal{R}_G{\upharpoonright} B$.
\end{theorem}

\emph{Sketch of proof.}
Since the action of $H$ is strongly ergodic, for every $\epsilon_0$, 
we may find $\delta _0 >0$ such that if $A\subseteq X$ is any Borel set satisfying $\sup _{s\in S} \mu (s^{-1}  A \triangle A ) < \delta _0$, then either $\mu (A)< \epsilon_0$ or $\mu (A) > 1-\epsilon_0$. 

Given $\epsilon>0$, we choose $\epsilon_0$ such that $\epsilon_0<\min\{\epsilon/8, 1/24\}$. 
Strong ergodicity for $H$ delivers $\delta_0$. We then choose $\delta$ satisfying the condition $\delta<\min\{\delta_0/2, 1-8\epsilon_0\}$.

By the commuting assumption, for every $k$ in the group $K$, for every $s$ in the generating set $S\subset H$ we check that $s^{-1}  A_k \triangle A_k \subseteq X\setminus (A_s\cap k^{-1}A_s )$.
Hence, by property (i),
$
\sup _{s\in S} \mu (s^{-1}  A_k \triangle A_k ) < 1-\mu (A_s\cap k^{-1}A_s )< 2\delta < \delta _0
$,
so that for each $k\in K$
\begin{equation}\label{eqn:orK}
\text{either \hskip10pt }\mu (A_k ) <\epsilon_0 \text{\hskip10pt  or \hskip10pt } \mu (A_k)>1-\epsilon_0.
\end{equation}
Consider now the subset $K_0 := \{ k\in K : \mu (A_k)> 1- \epsilon_0 \}$ of $K$.
\\-- Property (ii) along with \eqref{eqn:orK} and $\epsilon_0\leq \epsilon$,  imply $T\subseteq K_0$. 
\\
-- Since $\epsilon_0<1/3$, then $K_0$ is a subgroup of $K$. 
Indeed, clearly $K_0=K_0^{-1}$, and if $k_0,k_1\in K_0$ then $\mu (A_{k_0k_1})\geq \mu (A_{k_1}\cap k_1^{-1}A_{k_0}) > 1-2\epsilon_0>\epsilon_0$ hence $\mu (A_{k_0k_1})>1-\epsilon_0$ by \eqref{eqn:orK}, and thus $k_0k_1\in K_0$. 
 \\
 It follows that $K_0=K$.
We have shown that $\mu (A_k )>1-\epsilon_0$ for all $k\in K$.

 Theorem 2.7 of \cite{IKT09} then implies that $\mu (\{ x\in X : \psi x \, \mathcal{E} \, x \} ) > 1-4\epsilon_0$, for every element $\psi\in [\mathcal{R}_K]$ of the full group of the orbit equivalence relation $\mathcal{R}_K$ of $K$. Thus, by Lemma 2.14 of \cite{IKT09} there exists an $\mathcal{R}_K\cap \mathcal{E}$-invariant Borel set $B\subseteq X$ with $\mu (B)\geq 1-4\epsilon_0$ such that $\mathcal{R}_K{\upharpoonright} B\subseteq \mathcal{E}{\upharpoonright} B$. Indeed, $\mathcal{R}_K$ is \emph{relatively non-approximable in ${\mathcal{R}}_G$} (see below).
We now claim that 
\begin{equation}\label{claim:Gtime}
\textrm{for each $g\in G$, either $\mu (A_g) < 8\epsilon_0$, or $g^{-1}B\cap B\subseteq A_g$}\footnote{(thus in this case 
$\mu (A_g) >  1-8 \epsilon_0$)}
\end{equation}
If $\mu (A_g) > 8\epsilon_0$ for some $g\in G$. Then the set $A_g \cap g^{-1}B\cap B$ is a non-null subset of $B$, so it meets almost every $\mathcal{R}_K{\upharpoonright} B$ equivalence class since $\mathcal{R}_K{\upharpoonright} B$ is ergodic. For each $x\in g^{-1}B\cap B$ we can find some $k\in K$ such that $kx \in A_g\cap g^{-1}B\cap B$. Then $x,gx, kx,gkx \in B$ and $k,gkg^{-1}\in K$, so $x\, (\mathcal{R}_K{\upharpoonright} B) \, kx \, (\mathcal{E}{\upharpoonright} B) \, gkx =gkg^{-1}gx \, (\mathcal{R}_K{\upharpoonright} B) \, gx$, whence $x\in A_g$.

Let $G_0 = \{ g\in G : g^{-1}B\cap B \subseteq A_g \}$. 
\\
--
Since $8\epsilon_0<\epsilon$ and $1>1-\delta>8\epsilon_0$,
then properties (i) and (ii) and Claim (\ref{claim:Gtime})
 imply that $S\cup T \subseteq G_0$. 
\\
-- 
Since $\epsilon_0<1/24$ then $G_0$ is a subgroup of $G$: It is clear that $G_0^{-1}=G_0$ (since $A_{g^{-1}}=gA_g$). If $g_0,g_1 \in G_0$ then $\mu (A_{g_0})\geq 
1-8\epsilon_0$ and likewise $\mu (A_{g_1})\geq 1-8\epsilon_0$, so that $\mu (A_{g_0g_1})\geq \mu (A_{g_1}\cap g_1^{-1}A_{g_0})\geq 1-16\epsilon_0 >8\epsilon_0$ and hence $g_0g_1\in G_0$ by Claim~(\ref{claim:Gtime}).  
\\
Therefore, $G_0 = G$. This shows that $\mathcal{R}_G{\upharpoonright} B \subseteq \mathcal{E}{\upharpoonright} B$.
\qed

\medskip

Consider a pair ${\mathcal{S}}\subset {\mathcal{R}}$ of p.m.p. standard equivalence relations.
A standard sub-relation ${\mathcal{S}}\subset {\mathcal{R}}$ of p.m.p. standard equivalence relations is \emph{relatively non-approximable} if for every approximation $({\mathcal{R}}_n)$ of ${\mathcal{R}}$, there is some $n$ and a non-negligeable $A$ with ${\mathcal{S}}{\upharpoonright} A\subset {\mathcal{R}}_n{\upharpoonright} A$. This notion is useful through several variants of the following proposition.
\begin{proposition}[Weak form of normality]\label{prop: q-normality sub rel}
If ${\mathcal{R}}$ contains a sub-equivalence relation ${\mathcal{S}}$ and ${\mathcal{R}}$ is generated by a family $\phi_1, \phi_2, \cdots, \phi_p$ of isomorphisms of the space such that, $\phi_i({\mathcal{S}})\cap {\mathcal{S}}$ is ergodic for each~$i$.
Then every approximation $({\mathcal{R}}_n)$ for which there is a non-negligeable $A$ with ${\mathcal{S}}{\upharpoonright} A\subset {\mathcal{R}}_0{\upharpoonright} A$
has to be trivial. 
\end{proposition}
Consider such an approximation. We introduce the \emph{Window Trick}: 
\\
Let ${\mathcal{R}}'_n:= ({\mathcal{R}}_n{\upharpoonright} A)\vee {\mathcal{S}}$ be the sub-relation of ${\mathcal{R}}$ generated by ${\mathcal{R}}_n{\upharpoonright} A$ and ${\mathcal{S}}$. 
We claim that:
\\ (a) ${\mathcal{R}}'_n{\upharpoonright} A= {\mathcal{R}}_n{\upharpoonright} A$,
and 
\\ (b) $({\mathcal{R}}'_n)$ is an approximation of ${\mathcal{R}}$.

Now, the set $A_i^n:=\{x\in X: x{\mathcal{R}}'_n \phi_i^{-1}(x)\}$ is $(\phi_i({\mathcal{S}})\cap {\mathcal{S}})$-invariant: If $x\in A_i^n$ and $(x,y)\in \phi_i({\mathcal{S}})\cap {\mathcal{S}}$ then $y\overset{{\mathcal{S}}}{\sim}  x \overset{{\mathcal{R}}'_n}{\sim}\phi_i^{-1}(x) \overset{{\mathcal{S}}}{\sim} \phi_i^{-1}(y)$.
So that $y\in A_i^n$. It thus has full measure as soon as it is non-negligeable, and this happens for large enough $n$ since ${\mathcal{R}}'_n$ is an approximation.  
Taking an $n$ which is suitable for all $i$, we obtain ${\mathcal{R}}'_n={\mathcal{R}}$. 
So that ${\mathcal{R}}'_n{\upharpoonright} A= {\mathcal{R}}_n{\upharpoonright} A={\mathcal{R}}{\upharpoonright} A$. 
\qed

\medskip
Let $G=B(p,q)=\langle a, t \vert ta^pt^{-1}=a^{q}\rangle$ be a \emph{Bausmlag-Solitar group}.
The kernel $N$ of the modular map $G\to \mathbb{Q}^{*}, t\mapsto p/q, a\mapsto 1$ 
consists in the elements $w$ of $G$ which commute with a certain power $a^{k_w}$ of $a$.

\begin{theorem}[Baumslag-Solitar groups]\label{th: BS-groups}
If the kernel $N$ of the modular map acts strongly ergodically and all the (non trivial) powers of $a$ act ergodically, then the free action of $B(p,q)$ is non-approximable.
\end{theorem}

\medskip

Indeed, one can find a finitely generated subgroup $N_0$ of $N$ which already acts strongly ergodically. 
There is a common power $a^k$ which commutes with $N_0$. 
Applying Theorem~\ref{th:red herring fr} we show that $G_0=N_0. \langle a^k\rangle$ is non-approximable.
Thus the sub-relation generated by $G_0$ is relatively non-approximable. 
Proposition~\ref{prop: q-normality sub rel} applied to the pair of relations generated by $G_0$ and $G_1=N_0. \langle a\rangle$ with $\phi_1=a$, first; and then applied to the pair generated by
$G_1< B(p,q)$ with $\phi_1=t$ proves the result.\qed

We also obtain similar results for (most) inner amenable groups and various related families of groups.

\section{Approximate and geometric dimensions}
Besides consequences in Bernoulli bond percolation, Theorem~\ref{thm:first} allows us to obtain some information about the \emph{approximate dimension}.

A standard p.m.p. equivalence relation ${\mathcal{R}}$, when considered as a measured groupoid may act on bundles (fields) of simplicial complexes $x\mapsto \Sigma_x$ over $X$. The action is \emph{proper} if its restriction to the $0$-skeleton $x\mapsto \Sigma_x^{(0)}$ of the sub-bundle is smooth.
The \emph{dimension} of such a bundle is the maximum dimension of a fiber $\Sigma_x$, and the bundle is said to be \emph{contractible} if (almost) each fiber is contractible. 
The \emph{geometric dimension} of ${\mathcal{R}}$ is the minimum of the dimensions of the ${\mathcal{R}}$-bundles which are proper and contractible.
The \emph{approximate dimension} of ${\mathcal{R}}$ is the minimum of the dimensions $d$ such that ${\mathcal{R}}$ admits an approximation $({\mathcal{R}}_n)$ by sub-relations of dimension $d$.
These notions were introduced in \cite[section 5]{Gab02}.

For instance, smooth equivalence relations have geometric dimension $=0$.
Aperiodic treeable equivalence relations are exactly those with geometric dimension $=1$.
Their approximate dimension is $=0$ if and only if they are hyperfinite and is $=1$ otherwise.
One can show that the general inequalities:
approx-dim $\leq$ geom-dim $\leq$ approx-dim  $+1$.
It is unknown whether there are groups admitting free p.m.p. actions with different geometric dimensions. As for approximate dimension, various situations may occur. For instance, we obtain:
\begin{proposition}
Let $G_d:=\mathbf{F}_2\times \mathbf{F}_2\cdots\times \mathbf{F}_2\times \mathbb{Z}$ be the direct product of $d$ copies of the free group $\mathbf{F}_2$ and one copy of $\mathbb{Z}$. 
All its free p.m.p. actions have geometric dimension $=d+1$. It admits both free p.m.p. actions with approximate dimension $=d$ and $=d+1$.
\end{proposition}

\medskip

As already mentionned free products of equivalence relations are always approximable. This is 
no more the case for free actions of amalgamated free products over an infinite central subgroup $G=G_1*_C G_2$ when the common subgroup has indices greater than $3$ in the factors (apply Theorem~\ref{th:red herring fr} to, say, the Bernoulli shift action with $H=G$ and $K=C$).
This allows us to produce examples of group actions which are amalgamated free products of treeable over amenable, but which are not \emph{approxi-treeable}~(approximable by treeable): take for instance $G_1$ and $G_2$ abelian.

\section*{Acknowledgements}
The first author was supported by the ANR project GAMME (ANR-14-CE25-0004) and by the CNRS.
The second author was supported by NSF grant DMS 1303921.

\end{document}